\documentclass{amsart}

\usepackage{amsmath, amssymb, amsthm}
\usepackage{amsthm}
\usepackage{xcolor}
\usepackage{hyperref}
\usepackage{mathtools}

\newtheorem{thm}{Theorem}[section]
\newtheorem{prop}[thm]{Proposition}

\theoremstyle{definition}
\newtheorem{defn}[thm]{Definition}

\theoremstyle{remark}
\newtheorem{rem}[thm]{Remark}

\numberwithin{equation}{section}

\begin{document}

\title[A double-phase Neumann problem with $p=1$]{A double-phase Neumann problem with $p = 1$}
\author{Alexandros Matsoukas }
\address{Department of Mathematics\\
National Technical University of Athens\\
Iroon Polytexneiou 9\\
15780 Zografou\\
Greece}
\email{alexmatsoukas@mail.ntua.gr}
\author{Nikos Yannakakis}
\address{Department of Mathematics\\
National Technical University of Athens\\
Iroon Polytexneiou 9\\
15780 Zografou\\
Greece}
\email{nyian@math.ntua.gr}
\commby{}

\begin{abstract}
We study a double-phase Neumann problem with non-homogeneous boundary conditions, where the lowest exponent $p$ is equal to 1. The existence of a solution is established as the limit of solutions to corresponding double-phase problems with $p>1$. We also provide a variational characterization of the limit.
\end{abstract}

\maketitle

\section{Introduction}
\label{sec1}
In this paper we prove the existence of $W^{1,1}$-solutions for a double-phase Neumann problem involving the 1-Laplacian. In particular,  we study a boundary value problem with non-homogeneous Neumann condition, where the differential operator is the sum of the 1-Laplacian and a weighted q-Laplacian; that is a problem of the form
\begin{equation}\label{limiting}
\left\{
\begin{array}{rcl}
-\mathrm{div}\left(\frac{\nabla u}{|\nabla u|} +a(x)|\nabla u|^{q-2}\nabla u \right)&=&0\,\text{ in } \Omega,\\
\left[(\frac{\nabla u}{|\nabla u|} +a(x)|\nabla u|^{q-2}\nabla u),\nu \right]&=&g \, \text{ on }\partial\Omega,\,
\end{array}
\right.
\end{equation}
where $\Omega \subset \mathbb{R}^N$ is a bounded Lipschitz domain and $a(\cdot)$ is a bounded function with $a\ge 0$ a.e. in $\Omega$. As for the datum $g$ we assume that it belongs to $L^{\infty}(\partial \Omega)$ and satisfies the compatibility condition
\begin{equation*}
    \int_{\partial \Omega} g \, d \mathcal{H}^{N-1} =0.
\end{equation*}
Note that this latter condition is necessary for the well-posedness of problem \eqref{limiting}.  Moreover, since solutions to the Neumann problem are invariant under the addition of constants, we impose a normalization in order to recover uniqueness. To this end, we work in the subspace of the weighted Sobolev space consisting of functions with zero mean, namely
\begin{equation*}
    V(\Omega)=\bigg\{ u\in W_{a}^{1,q}(\Omega) : \int_{\Omega} u \, dx=0 \bigg\}.
\end{equation*}

The above problem, is a special case for $p = 1$, of $p,q$ double-phase problems i.e. problems where the differential operator is the sum of a $p$-Laplacian and a weighted $q$-Laplacian
\begin{equation}\label{pqdp}
\mathrm{div}\left(|\nabla u|^{p-2}\nabla u +a(x)|\nabla u|^{q-2}\nabla u \right)\,,\text{ for } \, u \in W^{1,\theta_p}(\Omega).\\
\end{equation}
This non-linear differential operator is related to the so-called double-phase functional
\begin{equation*}
 u \mapsto \int_{\Omega} (|\nabla u|^{p} +a(x)|\nabla u|^{q})\, d x, 
\end{equation*}
with $1<p<q$.

The double-phase functional was originally introduced by Zhikov \cite{Zhikov} in connection with the Lavrentiev phenomenon in the calculus of variations. These types of functionals were later investigated by Marcellini in the context of regularity theory for integral functionals with non-standard growth conditions of $p,q$ type \cite{Marcelini1}. The interest for studying double-phase problems emerged after the pioneering works of Mingione and co-workers \cite{mingione_1}, \cite{mingione_2}, where important local regularity results for minimizers of the double-phase functional were established. 
Since then, the terminology double-phase problems has become standard in the literature, as it effectively describes the characteristic feature of switching the growth of the functional in the two distinct regions of the domain determined by the weight function; that is on the sets 
$\{a>0\}$ and $\{a=0\}$.
Double-phase problems have attracted considerable attention and have been extensively studied by many authors; we refer the interested reader to the survey papers \cite{Chlebicka}, \cite{MinRad} \cite{P}, \cite{R} and the references therein. It is worth to mention that the majority of existing works focus on the superlinear growth case $p>1$, while the results for $p=1$ are quite sparse (see \cite{EHH1}, \cite{EHH2}, \cite{GLM}, \cite{harjulehto_1}, \cite{MY}, \cite{MY2}). We also mention here the works on perturbed 1-Laplacian equations by  p-Laplacian \cite{GTs,Ts}. The aim of the present work is to contribute further to the line of research initiated in \cite{MY,MY2}.

A relevant double-phase problem with homogeneous Neumann boundary conditions was previously addressed in \cite{GLM}, under a different set of assumptions and in a framework motivated by image restoration. The approach there, was based on convex duality methods and also includes the corresponding Neumann problem for the variable exponent case with $p(x) \ge 1$. For a comprehensive investigation of the effectiveness of the double-phase model with
$p=1$ in image restoration, together with a detailed exposition of numerical results, we refer to the recent work \cite{GLM2}. Additionally in \cite{KY}, a Neumann problem with non-homogeneous boundary conditions for the $p(x)$-Laplacian, where $p = 1$ in a subdomain, was studied. We refer the reader to the work \cite{MRST} on the Neumann problem involving the 1-Laplacian, for the non-homogeneous case, and to the foundational work \cite{ABCM} by Andreu, Ballester, Caselles, and Maz{\'o}n on the total variation flow. 
In particular, the study of problems for the pure 1-Laplacian and their approximation by 
$p$-Laplacian type equations as $p \to 1$ has a long history. Classical results in this direction include the works \cite{Mercaldo2008, Mercaldo2009} by Mercaldo, Segura de Le\'on, and Trombetti, which analyze the behavior of solutions to
$p$-Laplacian equations in the limit $p\to 1$ and address existence results for the 1-Laplacian with $L^{1}$-data, and more recently \cite{DeCicco2019, Latorre2021}, which treat Dirichlet problems involving singular nonlinearities. These works highlight the intrinsic singular nature of the 1-Laplacian, whose degeneracy is determined solely by the vanishing of the gradient. In contrast, the double-phase operator adds a second regularizing component driven by the weighted $q$-Laplacian. This interaction creates a competition between the linear-growth term and the superlinear $q$-growth term. Consequently, the structure of solutions and the associated variational framework can differ substantially from the pure 1-Laplacian case. In the context of double-phase problems involving the 1-Laplacian, boundary value problems with Dirichlet condition were studied in \cite{MY, MY2}, where the existence of solutions was obtained as the limit of solutions of approximate $p,q$ double-phase problems, as $p\to 1$.

In this paper we will follow the same approach in order to show that there exists a suitably defined weak solution to the Neumann problem \eqref{limiting}. We also give a variational characterization of this solution through the corresponding minimization problem.

 We first study the existence and behaviour of solutions for $p$ close to $1$ of the corresponding approximate double-phase Neumann problem with $p>1$ 
\begin{equation}\label{app}
\left\{
\begin{array}{rcl}
-\mathrm{div}\left(|\nabla u|^{p-2} \nabla u +a(x)|\nabla u|^{q-2}\nabla u \right)&=&0\,\text{ in } \Omega,\\
(|\nabla u|^{p-2}+a(x)|\nabla u|^{q-2}) \frac{\partial u}{\partial \nu} &=&g \, \text{ on }\partial\Omega.
\end{array}
\right.
\end{equation}
 Next, by considering a sequence of solutions $(u_p)$ to problems \eqref{app}, in order to show existence of solutions for the limiting problem \eqref{limiting} in the sense of a suitable definition, we pass to the limit as $p\to 1$. In particular, the solution $u$ to this problem is obtained as the limit, in the weak topology of the corresponding weighted Sobolev space, of the sequence of approximate solutions $(u_p)$, as $p \to 1$.

 As usual, in order to give sense to the quotient $\frac{\nabla u}{|\nabla u|}$, especially when the set $\{\nabla u = 0\}$ has positive Lebesgue measure, the above quotient will be replaced by a bounded vector field $z$ such that $\|z\|_{L^{\infty}}\le 1$ and which satisfies the condition $z \cdot \nabla u= |\nabla u|$ a.e. in $\Omega$. Another important issue to address is the interpretation of the boundary term corresponding to the Neumann condition. To this end, we introduce a notion of weak trace for the normal component of vector fields belonging to the class of double-phase compatible fields with integrable divergence
\begin{equation*}
X_{\mathrm {dp}}^N(\Omega) = \bigg\{ \zeta = z + a(x) |w|^{q-2} w:\, z\in L^\infty(\Omega)^N \,, w \in L_{a}^{q}(\Omega)^N \, |\, \mathrm{div}(\zeta) \in L^N(\Omega) \bigg\}.
\end{equation*}
More precisely, we define a function $\gamma (\zeta) \in W^{-\frac{1}{q},q}(\partial \Omega)$ which is associated to any vector field $\zeta \in X_{\mathrm {dp}}^{N}(\Omega )$ and we prove a weak Gauss-Green formula under the assumption that the trace of a function $u\in W^{1,q}_{a}(\Omega)$ belongs to $W^{1-\frac{1}{q},q}(\partial\Omega)$. This trace condition is ensured by requiring the weight function $a$ to be bounded away from zero on $\partial \Omega$.

Let us remark here that, as in the Dirichlet problem in \cite{MY}, we can study the approximate problems 
 \begin{equation*}
\left\{
\begin{array}{rcl}
-\mathrm{div}\left(|\nabla u|^{p-2} \nabla u +a(x)|\nabla u|^{q-2}\nabla u \right)&=&f\,\text{ in } \Omega,\\
(|\nabla u|^{p-2}+a(x)|\nabla u|^{q-2}) \frac{\partial u}{\partial \nu} &=&g \, \text{ on }\partial\Omega,
\end{array}
\right.
\end{equation*}
for $f\in L^{N}(\Omega)$ and $g\in L^{\infty}(\partial \Omega)$ and then pass to the limit as $p\to 1$. We point out that, our arguments still work just by considering a slightly different condition for the smallness of the norms of $f$ and $g$ in Proposition \ref{prop:behaviour} (see also \cite[Proposition 3.3]{MY}). For the sake of simplicity, we present here only the results for the case $f=0$.

\section{Notation and preliminary results}

In this section, we recall some definitions and tools of the function spaces we will use in our analysis.

\subsection{Weighted Lebesgue and Sobolev spaces}
We begin with the definition of the Muckenhoupt class $A_q$.
\begin{defn}
A weight $a\in L^\infty(\Omega)$ with $a(x)>0$ a.e. in $\Omega$ belongs to the Muckenhoupt class $A_q$ if
\[\sup_Q\left(\frac{1}{|Q|}\int_Q a(x) \, d x\right) \left(\frac{1}{|Q|}\int_Q a(x)^{-\frac{1}{q-1}} \, d x\right)^{q-1} < \infty,\]
where the supremum is taken over all cubes $Q$ with sides parallel to the coordinate axes.
\end{defn}
\noindent From now on we will always assume that the weight $a$ belongs to $A_q$. The weighted Lebesgue and Sobolev spaces are defined as
$$L^q_a (\Omega)=\{u: \Omega \to \mathbb{R}\, \text{measurable} \, : \int_{\Omega} a(x)|u|^q \, d x <+\infty\}$$
and
$$W^{1,q}_{a}(\Omega)=\left\{u\in L^{q}_{a}(\Omega): \, u \, \text{ is weakly differentiable and } |\nabla u|\in L^{q}_{a}(\Omega)\right\}\,.$$
Equipped with the norms 
$$\|u\|_{L^q_a} = \left( \int_{\Omega} a(x)|u|^q \, d x \right)^{\frac{1}{q}} \text{ and } \|u\|_{W^{1,q}_a}=\|u\|_{L^q_a}+\|\nabla u\|_{L^q_a}$$ $L^q_a (\Omega)$ and $W^{1,q}_{a}(\Omega)$  become uniformly convex and hence reflexive Banach spaces.

\noindent As usual, we define
\[W^{1,q}_{a, 0}(\Omega)=\overline{C_{0}^{\infty}(\Omega)}^{\|\cdot\|_{{W^{1,q}_a}}}\,.\]

We will also need the following Poincare-Wirtinger  inequality for the corresponding weighted Sobolev space with Muckenhoupt weight. 

\begin{prop}\label{prop:Poincare}
    Let $a\in A_q$ and $V(\Omega)\subset W^{1,q}_{a}(\Omega)$ be a closed subspace such that the only constant function is the zero function. Then we have 
\begin{equation}\label{PoincareAq}
    \|u\|_{L^q_a} \le c \|\nabla u\|_{L^q_a}, \,\,\, \text{ for all } \, u\in V(\Omega).
\end{equation}
\end{prop}

\begin{proof}
    The proof is by contradiction. Suppose that there exists a sequence  $(v_n)\in V_{a}(\Omega)$ such that 
    $$\|v_n\|_{L^q_a} > n \|\nabla v_n\|_{L^q_a}, \, \, \text{ for all } n\in \mathbb{N}.$$ Then, using the reflexivity of the space $W^{1,q}_{a}(\Omega)$, the compact embedding $W^{1,q}_{a}(\Omega) \hookrightarrow L^{q}_{a}(\Omega)$ and the weak lower semicontinuity of the norm $\|v\|_{W^{1,q}_{a}}$, we reach a contradiction.
\end{proof}
We will use the above Poincar{\'e}-Wirtinger inequality for the subspace
\begin{equation*}
    V(\Omega)=\bigg\{ u\in W_{a}^{1,q}(\Omega) : \int_{\Omega} u \, dx=0 \bigg\}.
\end{equation*}

If $\Omega$ is a bounded Lipschitz domain, $a\in C(\overline{\Omega})$ and is non-zero on $\partial\Omega$ then we can define a trace on $W^{1,q}_{a}(\Omega)$. 
\begin{prop}[\cite{MY}, Proposition 2.2]
\label{prop:trace}
Let $a\in C(\overline{\Omega})$ with $a\geq 0$ a.e. in $\Omega$, such that $a(x)\neq 0$, for all $x\in \partial\Omega$. Then there exists a bounded linear operator
\[T:W^{1,q}_{a}(\Omega)\rightarrow L^q(\partial\Omega)\]
such that
\[Tu=u|_{\partial \Omega} \text{ for all } u\in C(\overline{\Omega})\cap W^{1,q}_{a}(\Omega).\]
\end{prop}

Note that under the assumptions of Proposition \ref{prop:trace}, the trace of a function $v\in W^{1,q}_{a}(\Omega)$ belongs to the fractional Sobolev space $W^{1-\frac{1}{q},q}(\partial \Omega )$, which is defined as 
\begin{equation*}
    W^{1-\frac{1}{q},q}(\partial \Omega )= \bigg\{ v \in L^{q}(\partial \Omega) : \frac{|v(x)-v(y)|}{|x-y|^{(n-1)\frac{1}{q}+1}} \in L^{q}(\partial \Omega) \times L^{q}(\partial \Omega) \bigg\}.
\end{equation*}
It becomes a Banach space when endowed with the norm
$$ \|v\|_{W^{1-\frac{1}{q},q}(\partial \Omega )}= \|v\|_{L^{q}(\partial \Omega)}+ [v]_{W^{1-\frac{1}{q},q}(\partial \Omega )},$$
where $[v]_{W^{1-\frac{1}{q},q}(\partial \Omega )}$ denotes the so-called Gagliardo-Slobodeckij semi-norm 
\begin{equation*}
[v]_{W^{1-\frac{1}{q},q}(\partial \Omega )}= \bigg( \int_{\partial \Omega} \int_{\partial \Omega} \frac{|v(x)-v(y)|}{|x-y|^{n+q(1-\frac{1}{q})}} \, d s(x,y) \bigg)^{\frac{1}{q}}.
\end{equation*}
For more information in fractional Sobolev spaces, we refer the interested reader to \cite{DGP}.
\begin{rem}
\label{gagliardo}
   We remark here that, by a well-known theorem of Gagliardo (see \cite{Gagliardo}) for a Lipschitz domain $\Omega$, the trace operator is surjective from $W^{1,q}(\Omega)$ onto $W^{1-\frac{1}{q},q}(\partial \Omega)$. Moreover, there exists a bounded linear operator 
   $$\mathcal{E}:  W^{1-\frac{1}{q},q}(\partial \Omega )\to W^{1,q}(\Omega)$$ through which any function $\psi \in W^{1-\frac{1}{q},q}(\partial \Omega )$ can be extended to a function $v\in W^{1,q}(\Omega )$ such that $v|_{\partial \Omega}=\psi$, in the trace sense.
\end{rem}

The following Meyers-Serrin type approximation theorem, see \cite[Proposition 2.4]{MY2}, will play a significant role in what follows. 
  
\begin{prop}\label{prop:MS}
Assume that $a\in A_{q}$. If $u\in W^{1,1}(\Omega)\cap W^{1,q}_{a}(\Omega)$, then there exists a sequence $(v_n)$ in $W^{1,1}(\Omega)\cap C^\infty(\Omega)$ such that
\begin{align*}
v_n &\to u\,, \text{ in } W^{1,1}(\Omega),\\
\nabla v_n &\to \nabla u\,, \text{ in } L^{q}_{a}(\Omega)^N.
\end{align*}
If in addition $a(\cdot)$ is as in Proposition \ref{prop:trace} then we can actually take $(v_n)$ in $W^{1,q}(\Omega)\cap C^\infty(\Omega)$ and such that $v_{n}|_{\partial \Omega} = u|_{\partial \Omega}$, in the sense of the trace in $W^{1,q}_{a}(\Omega)$.
\end{prop}

\subsection{Generalized Orlicz spaces} For this part we follow mainly the survey paper \cite{P} and the book \cite{harjbook}. Let $\Omega \subset \mathbb{R}^N$ be a bounded Lipschitz domain, with $1<p<q<N$ and let $a\in L^{\infty}(\Omega)$  be a non-negative weight function. For fixed $q$ and with $p$ taking values in the above range, the functions
\[\theta_p :\Omega\times\mathbb R_+\rightarrow\mathbb R_+\]
defined by
$$ \theta_p(x,t)=t^p+a(x)t^q $$
are uniformly convex, generalized $\Phi$-functions \cite[Remark 2.22]{P} and satisfy the $(\Delta_2)$ condition \cite[Proposition 2.6]{P}.
The double-phase generalized Orlicz space is defined as
$$L^{\theta_p}(\Omega)=\{u: \Omega \to \mathbb{R}\, \text{measurable} : \rho_{\theta_p} (u) < +\infty \},$$
with modular given by
$$\rho_{\theta_p}(u)=\int_{\Omega} \theta_p(x,|u(x)|)\,d x.$$
When equipped with the so-called Luxemburg norm
$$\|u\|_{\theta_p} = \inf \bigg\{ \lambda>0 : \rho_{\theta_p} (\frac{u}{\lambda}) \le 1\bigg\},$$
$L^{\theta_p}(\Omega)$ becomes a uniformly convex (and hence reflexive) Banach space \cite[Proposition 2.23]{P}. 

\noindent The generalized Orlicz-Sobolev space is defined as
$$W^{1,\theta_p}(\Omega)=\left\{u\in L^{\theta_p}(\Omega): \,u \, \text{  is weakly differentiable and } \,|\nabla u|\in L^{\theta_p}(\Omega)\right\},$$
where $\nabla u$ is the weak gradient of $u$, and equipped with the norm
$$\|u\|_{1,\theta_p}=\|u\|_{\theta_p} + \|\nabla u\|_{\theta_p}\,.$$
is a reflexive Banach space.

\noindent As before, we define
\[W^{1,\theta_p}_0(\Omega)=\overline{C_{0}^{\infty}(\Omega)}^{\|\cdot\|_{1,\theta_p}}\,.\]

For the double-phase Orlicz space the following embeddings into the classical and weighted Lebesgue spaces hold

\begin{equation}
    L^{\theta_p}(\Omega) \hookrightarrow L^{p}(\Omega) \,\, \text{ and } \,\, L^{\theta_p}(\Omega) \hookrightarrow L_{a}^{q}(\Omega),
\end{equation}
and the constants of the embeddings are not greater than $1$.

 If $a\in C^{0,1}(\Omega)$ and $\frac{q}{p}<1+\frac{1}{N}$,
then the maximal operator is bounded on $L^{\theta_p}(\Omega)$ and the constant that bounds it is independent of $p$ (see \cite[Theorem 4.3.4]{harjbook}). Hence, by \cite[Theorem 6.2.8]{harjbook} we have the following Poincar{\'e}-Wirtinger type inequality 
\begin{equation}\label{Poincare}
    \|u-u_{\Omega}\|_{\theta_p} \le C_P \text{diam} (\Omega) \|\nabla u\|_{\theta_p}, \, \,\,\, \text{ for all } \, u \in W^{1,\theta_{p}}(\Omega),
\end{equation}
where $u_{\Omega}=\frac{1}{|\Omega|} \int_{\Omega} u \, dx $ and the constant $C_P$ depends only on the dimension and the John constant of the domain. 

For more details on generalized Orlicz spaces we refer the interested reader to the book \cite{harjbook}.

\section{Main Results}

Our assumptions on the weight function $\alpha(\cdot)$ and the exponents $1<p<q$ are the following.
\begin{equation*}
\mathrm{(H)}: a\in C^{0,1}(\overline{\Omega})\cap A_{q}, \, a(x) \neq 0 \text{ on } \partial \Omega, \, \text{ and } \, \frac{q}{p}<1+\frac{1}{N}.
\end{equation*}

\noindent Note that if $a\in A_q$, then $a^{-\frac{1}{q-1}}$ is locally integrable and hence $L^q_a (\Omega) \subseteq L^1_{loc} (\Omega)$. Thus if $\mathrm{(H)}$ is satisfied, then since near the boundary the weight is bounded away from zero we actually have that $W^{1,q}_{a} (\Omega) \hookrightarrow W^{1,1} (\Omega)$.

\subsection{Traces of double-phase vector fields with integrable divergence}

In \cite{Anzellotti} \cite{ChenFrid}, \cite{NovStra} and \cite{Temam}, under summability conditions on the divergence of a vector field, trace properties for its normal component on $\partial \Omega$ were obtained. In this section, we first introduce a space $X_{\mathrm {dp}}^{N}$ of compatible double-phase vector fields with integrable divergence. Then, we define a function $\gamma (\zeta) \in W^{-\frac{1}{q},q}(\partial \Omega)$ which is associated to any vector field $\zeta \in X_{\mathrm {dp}}^{N}(\Omega )$, where by $W^{-\frac{1}{q},q}(\partial \Omega)$ we denote the dual of $W^{1-\frac{1}{q},q}(\partial \Omega)$.
The set of double-phase compatible vector fields is 
\begin{equation*}
X_{\mathrm {dp}}^N(\Omega) = \bigg\{ \zeta = z + a(x) |w|^{q-2} w:\, z\in L^\infty(\Omega)^N \,, w \in L_{a}^{q}(\Omega)^N \, |\, \mathrm{div}(\zeta) \in L^N(\Omega) \bigg\}.
\end{equation*}
We also provide a weak Gauss-Green formula under the assumptions $\mathrm{(H)}$; in particular, this implies that the trace of a function $u\in W^{1,q}_{a}(\Omega)$ belongs to $W^{1-\frac{1}{q},q}(\partial\Omega)$.  Since the weight function $a$ is bounded we have that $X_{\text{dp}}^{N} \subset X_{q'}^{N}$ where $q'$ is as usual the conjugate exponent of $q$ and
$$X_{q'}^{N}(\Omega)=\{ \zeta \in L^{q'}(\Omega)^N\, |\, \mathrm{div} (\zeta) \in L^{N}(\Omega)\}.$$ Hence, the weak trace of the normal component can be defined for $\zeta \in  X_{\text{dp}}^{N}$ once it is defined for $\zeta \in X_{q'}^{N}$. 
When endowed with the norm
$$ \|\zeta\|_{q',N} = \|\zeta \|_{q'} + \|\mathrm{div} \zeta \|_N, $$ the space $(X_{q'}^{N},  \|\cdot\|_{q',N}) $ becomes a Banach space and $C^{\infty}_{0}(\overline{\Omega})^N $ is dense in $X_{q'}^{N}$ (see \cite{NovStra}). By Gagliardo's theorem, any function $\psi \in W^{1-\frac{1}{q},q}(\partial \Omega )$ can be extended to a function $v\in W^{1,q}(\Omega )$ such that $v|_{\partial \Omega}=\psi$, in the trace sense. We will apply this result in order to obtain an estimate on the bilinear form which is defined in the next theorem and consequently, this will enable us to define a weak trace of the normal component of a vector field $\zeta \in  X_{q'}^{N}$. We have the following

\begin{thm} There exists a bilinear map $\langle \zeta, v \rangle_{\partial \Omega}: X_{q'}^{N}(\Omega) \times W^{1,q}(\Omega) \to \mathbb{R}$ such that
\begin{align}\label{bilestimate}
|\langle \zeta, v \rangle_{\partial \Omega}| \le  C' (  \|\zeta\|_{q'} +  \|\mathrm{div} (\zeta)\|_{N})  \|\psi \|_{1-\frac{1}{q},q}, \text{ for all } \zeta , v,
\end{align}
where $\psi \in W^{1-\frac{1}{q},q}(\partial \Omega)$ is such that $v|_{\partial \Omega}=\psi $.
Any vector field $\zeta \in X_{q'}^{N}(\Omega)$ defines a linear functional $$F_{\zeta} :  W^{1-\frac{1}{q},q}(\partial \Omega) \to \mathbb{R}$$ through the formula 
\begin{equation}\label{defbi}
\nonumber
F_{\zeta}(\psi)=\langle \zeta, v \rangle_{\partial \Omega}.
\end{equation}
 In particular, we can define a weak trace $\gamma(\zeta) \in W^{-\frac{1}{q},q}(\partial \Omega)$, of the normal component of $\zeta$ on $\partial \Omega$, such that the following integration by parts formula holds 
\begin{eqnarray}\label{ibp1}
\int_{\Omega} \zeta \cdot \nabla v \, d x+\int_{\Omega} v\, \mathrm{div} (\zeta)  \, d x =  \langle  \gamma (\zeta), \psi \rangle_{W^{-\frac{1}{q},q} \times W^{1-\frac{1}{q},q}}
\end{eqnarray}
for all $\zeta \in X_{q'}^{N}$ and $v\in W^{1,q}(\Omega)$, with $\psi=v |_{\partial \Omega}$.
\end{thm}

\begin{proof} Let $\psi \in W^{1-\frac{1}{q},q}(\partial \Omega)$ and $v\in W^{1,q}(\Omega)$ be such that $v|_{\partial \Omega}=\psi $. For $\zeta \in  X_{q'}^{N}(\Omega)$ we set 
\begin{equation}\label{bilinear}
\nonumber
\langle \zeta, v \rangle_{\partial \Omega} =\int_{\Omega} \zeta \cdot \nabla v \, d x + \int_{\Omega} v \, \mathrm{div}(\zeta)  \, d x.
\end{equation}
The map $\langle \zeta, v \rangle_{\partial \Omega}$ is bilinear and we observe that if $v_1,v_2 \in W^{1,q}(\Omega)$ and $v_{1}|_{\partial \Omega}=v_{2}|_{\partial \Omega}=\psi $, then 
\begin{eqnarray*}
\langle \zeta, v_1 \rangle_{\partial \Omega}=\langle \zeta, v_2 \rangle_{\partial \Omega} \text{ for all } \zeta \in X_{q'}^{N}(\Omega).
\end{eqnarray*}
To see this note that if $v=v_{1}-v_{2}$, then $v \in W_{0}^{1,q}(\Omega)$ and thus there exists a sequence $v_{n} \in C^{\infty}_{0}(\Omega)$ converging to $v$. So by the definition of the distributional divergence we get
\begin{align*}
\int_{\Omega} \zeta \cdot \nabla (v_{1}-v_{2}) \, d x &+ \int_{\Omega} (v_{1}-v_{2}) \mathrm{div}( \zeta) \,  d x \\ &=  \lim_{n} \bigg(\int_{\Omega} \zeta \cdot \nabla  v_{n} \, d x + \int_{\Omega} v_{n} \mathrm{div}(\zeta) \, d x \bigg)\\
&= 0\,.
\end{align*}

For $v\in W^{1,q}(\Omega)$ and $\zeta \in  X_{q'}^{N}(\Omega)$  by H\"older's inequality, and the Sobolev embedding theorem $W^{1,1}(\Omega) \hookrightarrow L^{\frac{N}{N-1}}(\Omega)$ we get
\begin{eqnarray*}
|\langle \zeta, v \rangle_{\partial \Omega}| &=&|\int_{\Omega} \zeta \cdot \nabla v \, d x + \int_{\Omega}v \,  \mathrm{div} (\zeta ) dx | \\
&\le & \|\zeta\|_{q'} \|\nabla v\|_{q} + \|\mathrm{div} (\zeta)\|_{N} \|v\|_{\frac{N}{N-1}} \\
&\le & \|\zeta\|_{q'} \|\nabla v\|_{q} + \|\mathrm{div} (\zeta)\|_{N}  |\Omega |^{1-{\frac{1}{q}}} C_{S} \|v\|_{1,q} \\
&\le & C (  \|\zeta\|_{q'} +  \|\mathrm{div} (\zeta)\|_{N})  \|v\|_{1,q}, \,
\end{eqnarray*}
where $C=\max\{1, |\Omega |^{1-{\frac{1}{q}}} C_{S}\}$ and $C_{S}$ is the constant of the Sobolev inequality. By Gagliardo's theorem there exists $c>0$ such that $\|v\|_{1,q} \le c \|\psi \|_{1-\frac{1}{q},q}$ whenever $v\in W^{1,q}(\Omega)$ is such that $v|_{\partial \Omega}=\psi$. Combining this with the above we get for such $v$ and $\psi$ that
\begin{eqnarray*}
|\langle \zeta, v \rangle_{\partial \Omega}|
\le  C' (  \|\zeta\|_{q'} +  \|\mathrm{div} (\zeta)\|_{N})  \|\psi \|_{1-\frac{1}{q},q}.
\end{eqnarray*}
Now, for $\zeta \in  \mathcal{X}_{q'}^{N}(\Omega)$  we define the linear functional $F_{\zeta} :  W^{1-\frac{1}{q},q}(\partial \Omega) \to \mathbb{R}$ with
\begin{equation*}
F_{\zeta}(\psi)=\langle \zeta, v \rangle_{\partial \Omega},
\end{equation*}
where $\psi \in W^{1-\frac{1}{q},q}(\partial \Omega)$ and $v\in W^{1,q}(\Omega)$ are such that $v|_{\partial \Omega}=\psi $. By the estimate \eqref{bilestimate} we have that 
\begin{equation*}
|F_{\zeta}(\psi)|\le C' (  \|\zeta\|_{q'} +  \|\mathrm{div} (\zeta)\|_{N})  \|\psi \|_{1-\frac{1}{q},q} 
\end{equation*}
and hence, there exists a function $\gamma(\zeta) \in  W^{-\frac{1}{q},q}(\partial \Omega)$ such that
\begin{equation*}
F_{\zeta}(\psi)= \langle  \gamma (\zeta), \psi \rangle_{W^{-\frac{1}{q},q} \times W^{1-\frac{1}{q},q}}
\end{equation*}
with
\begin{equation*}
\|\gamma (\zeta)\|_{W^{-\frac{1}{q},q}} \le c \|\zeta\|_{q',N}.
\end{equation*}
As a consequence, the integration by parts formula \eqref{ibp1} is immediate from the way we defined the functional $F_{\zeta}$.
\end{proof}

Eventually, we can use the integration by parts formula \eqref{ibp1}, in order to obtain the following Gauss-Green formula for double-phase vector fields with integrable divergence.

\begin{prop}
   Assume that $\mathrm{(H)}$ holds and let $u\in W^{1,q}_{a}(\Omega)$ and $\zeta \in X_{\rm dp}^{N}(\Omega ) $. Then one has the Gauss-Green formula, 
\begin{equation}\label{GaussGreen}
\int_{\Omega } u \,\mathrm{div}(\zeta) \,  dx + \int_{\Omega } z \cdot \nabla u \, dx + \int_{\Omega } a(x) |w|^{q-2}w \cdot \nabla u \, dx = \langle  \gamma (\zeta),  u \rangle_{W^{-\frac{1}{q},q} \times W^{1-\frac{1}{q},q}},
\end{equation}
where the duality pairing between the spaces $W^{-\frac{1}{q},q}(\partial \Omega)$ and $W^{1-\frac{1}{q},q}(\partial \Omega)$, lacks an integral representation.
\end{prop}

\begin{proof} Take a sequence of functions $v_n \in C^\infty (\Omega) \cap W^{1,q} (\Omega)$ that converges to $u$ as in Theorem \ref{prop:MS}. By the integration by parts formula (\ref{ibp1}), for $v_{n}$ in the place of $v$, we have that
\begin{eqnarray*}
\int_{\Omega }  v_n \, \mathrm{div} (\zeta) \, d x + \int_{\Omega } \zeta \cdot \nabla v_n \, d x =  \langle  \gamma (\zeta),  v_n \rangle_{W^{-\frac{1}{q},q} \times W^{1-\frac{1}{q},q}}\,,
\end{eqnarray*}
for all $n\in\mathbb N$.
Since $W^{1,1}(\Omega ) \hookrightarrow L^{\frac{N}{N-1}}(\Omega )$ and $\mathrm{div}( \zeta) \in L^{N}(\Omega )$ we get
\begin{eqnarray*}
\lim_n \int_{\Omega }  v_n \, \mathrm{div} (\zeta) \, d x = \int_{\Omega }  u \, \mathrm{div}(\zeta) \, d x.
\end{eqnarray*}
We know that, $\nabla v_n \to \nabla u $ in $L^{q}_{a}(\Omega)^N$ and hence $a^\frac{1}{q} \nabla v_n \to a^\frac{1}{q} \nabla u $ in $L^{q}(\Omega)^N$. Since, $a^\frac{1}{q'} |w|^{q-2}w \in L^{q'}(\Omega)^N$ we can now pass to the limit to get
\begin{eqnarray*}
\lim_n \int_{\Omega } \zeta \cdot \nabla v_n \, d x &=& \lim_n \int_{\Omega } z \cdot \nabla v_n \, d x + \lim_n \int_{\Omega } a(x) \, |w|^{q-2}w \cdot \nabla v_n \, d x\\
&=& \int_{\Omega } z \cdot \nabla u \, d x + \int_{\Omega } a(x) \, |w|^{q-2}w \cdot \nabla u \, d x\,,
\end{eqnarray*}
because $z \in L^{\infty}(\Omega)^N$ and $\nabla v_n \to \nabla u$ in $L^1 (\Omega)^N$. 
Eventually since, $v_n |_{\partial \Omega }=u|_{\partial \Omega }$, we have
\begin{eqnarray*}
\lim_n  \langle  \gamma (\zeta), v_n \rangle_{W^{-\frac{1}{q},q} \times W^{1-\frac{1}{q},q}} = \langle  \gamma (\zeta), u \rangle_{W^{-\frac{1}{q},q} \times W^{1-\frac{1}{q},q}}  
\end{eqnarray*}
and the proof is completed.
\end{proof}

\subsection{Existence of solutions to the limiting problem}

In this section, we will show the existence of solutions to the limiting double-phase Neumann problem involving the 1-Laplacian. This will be achieved by studying the behaviour of the weak solutions $(u_p)$ as $p \to 1$ of the approximate Neumann problems
\begin{equation}\label{approximate}
\left\{
\begin{array}{rcl}
-\mathrm{div}\left(|\nabla u|^{p-2} \nabla u +a(x)|\nabla u|^{q-2}\nabla u \right)&=&0\,\text{ in } \Omega,\\
(|\nabla u|^{p-2}+a(x)|\nabla u|^{q-2}) \frac{\partial u}{\partial \nu} &=&g \, \text{ on }\partial\Omega,\,
\end{array}
\right.
\end{equation}
As usual, we say that $u\in W^{1,\theta_{p}}(\Omega)$ is a weak solution to problem \eqref{approximate}, if 
\begin{equation}\label{weakform}
\int_{\Omega} \left(|\nabla u|^{p-2} \nabla u +a(x)|\nabla u|^{q-2}\nabla u \right) \nabla v \, dx = \int_{\partial \Omega} g v \, d\mathcal{H}^{N-1} \,, 
\end{equation}
for all $ v\in W^{1,\theta_{p}}(\Omega)$.

In the next Proposition we prove existence of solutions to the approximate double-phase Neumann problems using the direct method of the calculus of variations. In what follows, $\Lambda_{1}(\Omega)$ is the best constant of the trace embedding $W^{1,1}(\Omega) \hookrightarrow L^{1}(\partial \Omega)$ (see \cite{AMR}). 
 
\begin{prop}\label{existencepq}
Let $g\in L^{\infty}(\partial \Omega)$ and assume that $\mathrm{(H)}$ holds. Then, there exists a unique weak solution $u \in W^{1,\theta_p}(\Omega)$ to the double-phase Neumann problem \eqref{approximate}, which is the unique minimizer of the functional 
\begin{equation*}
   \mathcal{F}(u) = \int_{\Omega} (\frac{|\nabla u|^{p}}{p} +a(x)\frac{|\nabla u|^{q}}{q})\, d x - \int_{\partial \Omega} g u \, d\mathcal{H}^{N-1},
\end{equation*}
    in the space
\begin{equation*}    
    V_p (\Omega)= \{ u\in W^{1,\theta_p}(\Omega) : \int_{\Omega}u \, dx=0\}.
\end{equation*}   
\end{prop}

\begin{proof}
   In order use the direct method of the calculus of variations, we need to show that the functional $\mathcal{F}$ is coercive and sequentially weakly lower semicontinuous in $V_p(\Omega)$. By Poincare's inequality \ref{Poincare}, $\|\nabla u\|_{L^{\theta_p}}$ is an equivalent norm in $V_p(\Omega)$. For coercivity, we will show that if $\|\nabla u\|_{L^{\theta_p}} \to +\infty$ then $\mathcal{F}(u)\to +\infty$ hence, we may assume that $\|\nabla u\|_{L^{\theta_p}}>1$. We estimate $\mathcal{F}(u)$ from below, using the embedding $W^{1,\theta_p}(\Omega) \hookrightarrow W^{1-\frac{1}{q},q}(\partial \Omega)$, Poincare's inequality \eqref{Poincare} to get that
   \begin{eqnarray*}
       \mathcal{F}(u) &=& \int_{\Omega} (\frac{|\nabla u|^{p}}{p} +a(x)\frac{|\nabla u|^{q}}{q})\, d x - \int_{\partial \Omega} gu \, d\mathcal{H}^{N-1} \\
       &\ge& \frac{1}{q} \rho_{\theta_p} (|\nabla u|) - \|g\|_{\infty} \|u\|_{L^1(\partial \Omega)} \\
       &\ge & \frac{1}{q} \rho_{\theta_p} (|\nabla u|) - \Lambda_1(\Omega) |\Omega|^{1-\frac{1}{p}} \|g\|_{\infty}\|u\|_{W^{1,p}(\Omega)} \\
       &\ge & \frac{1}{q} \rho_{\theta_p} (|\nabla u|) - \Lambda_1(\Omega) |\Omega|^{1-\frac{1}{p}} \|g\|_{\infty}\|u\|_{1,\theta_p} \\
       &\ge & \frac{1}{q \alpha}  \|\nabla u\|^p_{\theta_p} - \Lambda_1(\Omega) |\Omega|^{1-\frac{1}{p}}(C_P \text{diam}(\Omega)+1) \|g\|_{\infty} \|\nabla u\|_{\theta_p} \,
   \end{eqnarray*}
   where we passed from the modular to the norm using \cite[Lemma 3.2.9]{harjbook} for some suitable constant $\alpha$. Since $p>1$, we have that $\mathcal{F}(u)\to +\infty $ as $\|\nabla u\|_{L^{\theta_p}} \to +\infty$, i.e. $\mathcal{F}$ is coercive. Next, we will show that $\mathcal{F}$ is sequentially weakly lower semicontinuous. If $u_n \stackrel{w}{\rightarrow} u$ in $V_p(\Omega)$, then by the weak to weak continuity of the trace 
   \begin{equation*}
   \int_{\partial \Omega} g u_n \, d\mathcal{H}^{N-1} \to \int_{\partial \Omega} g u \, d\mathcal{H}^{N-1}.
   \end{equation*}
   By weak lower semicontinuity \cite[Theorem 2.2.8]{DHHR} we get that 
\begin{equation*}
\int_{\Omega} (\frac{|\nabla u|^{p}}{p} +a(x)\frac{|\nabla u|^{q}}{q})\, d x  \le \liminf_{n}{ \int_{\Omega} (\frac{|\nabla u_n|^{p}}{p}+a(x)\frac{|\nabla u_n|^{q}}{q})\,d x  }
\end{equation*}
and hence we conclude that 
\begin{equation*}
\mathcal{F}(u) \le \liminf_{n}{\mathcal{F}}(u_n). 
\end{equation*}
Hence, the direct method of the calculus of variations produces a minimizer of the functional in the space $V_p(\Omega)$, which is unique due to strict convexity.
\end{proof}

From now on, with a slight abuse of notation, we will say that $(u_p)$ is a sequence and we will consider subsequences of it, for $p$ close to $1$. In the next proposition, we will study the behaviour of the solutions $(u_p)$ to the approximate problems \eqref{approximate}, in order to pass to the limit as $p\to 1$.

\begin{prop}[Behaviour of solutions as p goes to 1]\label{prop:behaviour}
    Let $g\in L^{\infty}(\partial \Omega)$, assume that $\mathrm{(H)}$ holds and 
    \begin{eqnarray}\label{datanorm}
        2\Lambda _{1}(\Omega)  (C_P \mathrm{diam}(\Omega) +1) \|g\|_{\infty} <1.
    \end{eqnarray}
    If $(u_p)$ are the weak solutions of problems \eqref{approximate}, then there exist a function $u\in  W^{1,q}_{a}(\Omega)$ and a vector field $z \in L^{\infty}(\Omega)^N$, with $\|z\|_{\infty}\le 1$ such that as $p\rightarrow 1$, up to subsequences 
\begin{align*}
    u_p &\to u \, \text{ in } L^{s}(\Omega), \text{ for all } 1\le s<\frac{N}{N-1}\,,\\
    |\nabla u_p |^{p-2} \nabla u_p &\stackrel{w}{\rightarrow} z \, \text{ in } L^{r}(\Omega)^N \,, \text{ for all } \, 1\le r < +\infty,\\
 |\nabla u_p|^{q-2} \nabla u_p &\stackrel{w}{\rightarrow}  |\nabla u|^{q-2} \nabla u \, \text{ in } \, L^{q'}_{a}(\Omega)^N\,,\\
 \nabla u_p &\to \nabla u \, \text{ in } \, L^{q}_{a}(\Omega)^N\,.\\
\end{align*}

\end{prop}

\begin{proof}
    Let $\lambda_p = \|\nabla u_p \|_{L^{\theta_p}}$ and take $\frac{u_p}{\lambda_p}$ as a test function in the weak formulation (\ref{weakform}). Then, using the continuity of the trace operator 
     $W^{1,\theta_p}(\Omega) \hookrightarrow W^{1,1} (\Omega) \hookrightarrow L^{1}(\partial \Omega)$ and the Poincar{\'e} inequality (\ref{Poincare}), we have the following estimates
\begin{eqnarray}
  \nonumber  \frac{1}{\lambda_p} \int_{\Omega} (|\nabla u_p|^p +a(x)|\nabla u_p|^q) \, d x &=& \frac{1}{\lambda_p} \int_{\partial \Omega} g u_p \, d\mathcal{H}^{N-1} \\
 \nonumber   &\le &  \frac{1}{\lambda_p}  \|g\|_{\infty} \|u_p\|_{L^{1}(\partial \Omega)}\\
 \nonumber    &\le &  \frac{1}{\lambda_p} \Lambda_1 (\Omega) \|g\|_{\infty} \|u_p\|_{W^{1,1}(\Omega)} \\
 \nonumber    &\le &  \frac{1}{\lambda_p} \Lambda_1 (\Omega) |\Omega|^{1-\frac{1}{p}}\|g\|_{\infty} \|u_p\|_{1,p} \\
 \nonumber    &\le &  \frac{2}{\lambda_p} \Lambda_1 (\Omega) \|g\|_{\infty} \|u_p\|_{1,\theta_p} \\
 \nonumber   &\le & \frac{2}{\lambda_p} \Lambda_1 (\Omega) \|g\|_{\infty} (C_P \text{diam}(\Omega) +1) \|\nabla u_p\|_{\theta_p} \\
  \label{rhs}  &= & 2\Lambda_1 (\Omega) \|g\|_{\infty} (C_P \text{diam}(\Omega) +1) \,,
\end{eqnarray}
for $p$ close to 1. Now, we estimate also the left hand side 
\begin{eqnarray}
\nonumber  \frac{1}{\lambda_p} \int_{\Omega} (|\nabla u_p|^p +a(x)|\nabla u_p|^q) \, d x
\end{eqnarray}
\begin{eqnarray}    
\nonumber &=& {\lambda_p}^{p-1} \int_\Omega \bigg|\frac{\nabla u_p}{\lambda_p}\bigg|^pd x + {\lambda_p}^{q-1} \int_{\Omega} a(x) \bigg|\frac{\nabla u_p}{\lambda_p}\bigg|^q\, d x \\
\nonumber  &\ge &   \min\{\lambda_p^{p-1},\lambda_p^{q-1}\} \int_\Omega \bigg(\bigg|\frac{\nabla u_p}{\lambda_p}\bigg|^p + a(x) \bigg|\frac{\nabla u_p}{\lambda_p}\bigg|^q\bigg)\, d x \\
\label{lhs} &= &\min\{\lambda_p^{p-1},\lambda_p^{q-1}\}.\,
\end{eqnarray}
Thus, by \eqref{rhs} and \eqref{lhs}, for $p$ close to 1, we get that
\begin{equation*}
     \min\{\lambda_p^{p-1},\lambda_p^{q-1}\}\le 2\Lambda_1 (\Omega) \|g\|_{\infty} (C_P \text{diam}(\Omega) +1)
\end{equation*}
and by the assumptions on the norm of the data \eqref{datanorm}, we have that
\begin{equation*}
    \|\nabla u_p \|_{\theta_p} < 1.
\end{equation*}
Using the embedding $L^{\theta_p}(\Omega) \hookrightarrow L^{q}_a(\Omega)$ we also have
\begin{equation}
\nonumber
\left( \int_\Omega a(x) |\nabla u_p|^q \right)^{\frac{1}{q}} \le \|\nabla u_p\|_{\theta_p} < 1
\end{equation}
and hence by Poincar{\'e}'s inequality \eqref{PoincareAq} we have that $(u_p)$ is  bounded in $W^{1,q}_{a} (\Omega)$. By reflexivity there exists $u\text{ in } W_{a}^{1,q}(\Omega)$ such that, after passing to a subsequence, $$u_p \stackrel{w}{\rightarrow} u\,,\text{ in } W_{a}^{1,q}(\Omega)\,.$$
As we have already mentioned since $\mathrm{(H)}$ holds we have that $W^{1,q}_{a} (\Omega) \hookrightarrow W^{1,1} (\Omega)$ and hence $u_p \stackrel{w}{\rightarrow} u$ in $W^{1,1} (\Omega)$.  
Thus since for $1\le s <\frac{N}{N-1}$ the embedding $W^{1,1} (\Omega)\hookrightarrow L^{s}(\Omega)$ is compact, we have that
\begin{equation*}
 u_p\rightarrow u \, \text{ in } L^s(\Omega)\,,\text{ for all }1\le s <\frac{N}{N-1}\,.
\end{equation*} 
  Let $1\le r < p' $. Then by $L^{\theta_p}(\Omega) \hookrightarrow L^p(\Omega)$ and the fact that $\|\nabla u_p\|_{L^{\theta_p}} <1$, for $p$ close to 1, we obtain
\begin{eqnarray}
\nonumber \int_\Omega |\nabla u_p|^{(p-1)r}\,dx &\le &|\Omega|^{1-\frac{(p-1)r}{p}}\left( \int_\Omega |\nabla u_p|^{p}\,dx \right)^{\frac{(p-1)r}{p}}\\
\nonumber &\le & |\Omega|^{1-\frac{(p-1)r}{p}} \|\nabla u_p\|_{\theta_p}^{(p-1)r}\\
\nonumber &\le & |\Omega|^{1-\frac{(p-1)r}{p}}\\
\label{infinity} &\le &(|\Omega|+1)\,.
\end{eqnarray}
Note that for any fixed $r\geq 1$ by taking $p$ close enough to $1$ we get that $1\leq r<p'$. Hence by \eqref{infinity} the sequence  $(|\nabla u_p|^{p-2} \nabla u_p)$ is bounded in $L^r(\Omega)^N$ and thus it converges weakly to a $z_r\in L^r(\Omega)^N$. By a diagonal argument we may find a subsequence and a common vector field $z\in L^r(\Omega)^N$ such that
$$|\nabla u_p|^{p-2}\nabla u_p\stackrel{w}{\rightarrow} z\, \text{ in }L^r(\Omega)^N\,,$$
for all $1\le r < +\infty$. By \eqref{infinity} and using the fact that the norm is weakly lower semicontinuous we get that
\[\|z\|_r\leq (1+|\Omega|)^{\frac{1}{r}}\]
and thus
\[\|z\|_\infty=\lim_{r\rightarrow \infty}\|z\|_r\leq 1\,.\]
The last two convergences of the Proposition follow with a slight modification of the argument in the proof of \cite[Proposition 3.6]{MY}.  To this end let $\varepsilon>0$. Then by Proposition \ref{prop:MS} there exists $v\in W^{1,q}(\Omega)\cap C^\infty(\Omega)$ with $u|_{\partial \Omega}=v|_{\partial \Omega}$ such that
$$\bigg| \int_{\Omega} |\nabla v|\,dx-\int_{\Omega} |\nabla u| \, dx\bigg|<\frac{\varepsilon}{2}$$
and
\begin{equation}
\label{q_smooth_estimate}
(\int_\Omega a(x)|\nabla u-\nabla v|^q\, dx)^\frac{1}{q}<\frac{\varepsilon}{2}\,.
\end{equation}
Since $W^{1,q}(\Omega)\hookrightarrow W^{1,\theta_p}(\Omega)$ we may use $u_p-v$ in the weak formulation \eqref{weakform} and obtain
$$\int_\Omega (|\nabla u_p|^{p-2}\nabla u_p+a(x)|\nabla u_p|^{q-2}\nabla u_p) \nabla (u_p-v)\, dx=\int_{\partial \Omega}g (u_p-v) \, d\mathcal{H}^{N-1},$$
or equivalently,
\[\int_\Omega |\nabla u_p|^p\,dx-\int_\Omega |\nabla u_p|^{p-2}\nabla u_p\nabla v\,dx+\int_\Omega a(x)|\nabla u_p|^{q-2}\nabla u_p \nabla (u_p-v) \, dx=\]
$$\int_{\partial \Omega}g (u_p-v) \, d\mathcal{H}^{N-1}.$$
By Young's inequality we obtain
\[
\int_\Omega|\nabla u_p| \le \frac{1}{p}\int_{\Omega}|\nabla u_p|^p +\frac{p-1}p|\Omega|\,,
\]
and so we have
\[
p\int_\Omega |\nabla u_p|\,dx-\int_\Omega |\nabla u_p|^{p-2}\nabla u_p\nabla v\,dx+\int_\Omega a(x)|\nabla u_p|^{q-2}\nabla u_p \nabla (u_p-v) \,dx=\]
\[=\int_{\partial \Omega}g (u_p-v) \, d\mathcal{H}^{N-1}+(p-1)|\Omega|\,.\]
Letting $p\to 1$, using the lower semicontinuity of the norm, the fact that $|\nabla u_p|^{p-2}\nabla u_p \overset{w}{\rightharpoonup} z$ and that $u_p \overset{w}{\rightharpoonup} u$ in $W^{1,q}_{a}(\Omega)$ implies that the traces $Tu_p \overset{w}{\rightharpoonup} Tu$, we get that
\[\int_{\Omega}|\nabla u| \, dx+\limsup_{p\to 1}\int_\Omega a(x)|\nabla u_p|^{q-2}\nabla u_p \nabla (u_p-v) \,dx\le\]
\[\le\int_{\partial \Omega}g (u-v) \, d\mathcal{H}^{N-1}+\int_\Omega z \cdot \nabla v\,dx=\int_\Omega z \cdot \nabla v\,dx\,,\]
since $u|_{\partial \Omega}=v|_{\partial \Omega}$. Now since $\|z\|_{\infty}\le 1$ we get
\[
\bigg|\int_\Omega z \cdot \nabla v\, dx \bigg|\le \|z\|_{L^\infty} \int_\Omega |\nabla v|\, dx\le \int_\Omega |\nabla v|\, dx\,\]
and thus we get that
\begin{eqnarray}
\nonumber
\limsup_{p\to 1}\int_\Omega a(x)|\nabla u_p|^{q-2}\nabla u_p \nabla (u_p-v) \, dx
\le \int_\Omega|\nabla v|\, dx-\int_\Omega |\nabla u|\, dx<\frac{\varepsilon}{2}\,.
\end{eqnarray}
The rest follows the lines of the proof of \cite[Proposition 3.6]{MY} using the fact that $L^q_a(\Omega)$ being uniformly convex has the Radon-Riesz property.
\end{proof}

Now, we give a suitable definition of a weak solution to problem \eqref{limiting}.

\begin{defn}\label{defsol}
    We say that $u\in  W^{1,q}_{a}(\Omega)$ is a solution to problem \eqref{limiting} if there exists $z\in L^{\infty}(\Omega )^N$ with $\|z\|_{\infty}\le 1$ such that
    \begin{eqnarray*}
\int_{\Omega} z \cdot \nabla \phi \, d x+ \int_{\Omega} a(x) |\nabla u|^{q-2} \nabla u \cdot \nabla \phi \, d x &=& 0\,,\text{ for all }\phi\in C^\infty_0(\Omega),\\
z \cdot \nabla u &=& |\nabla u| \, \text{ a.e. in } \Omega,\\
\gamma(\zeta)&=&g \,\,\,\, \mathcal{H}^{N-1} \, \text{- a.e. on } \, \partial \Omega, 
\end{eqnarray*}
where $\zeta = z +a(x) |\nabla u|^{q-2}\nabla u$. \
\end{defn}

\begin{rem}[Weak formulation]\label{weaklimiting}
    If $u\in W^{1,q}_{a}(\Omega)$ is a weak solution to problem \eqref{limiting}, then it satisfies
    \begin{equation}\label{weakformlim}
        \int_{\Omega} |\nabla u| \, d x-\int_{\Omega} z \cdot \nabla v \, d x + \int_{\Omega} a(x) |\nabla u|^{q-2} \nabla u \cdot \nabla (u-v) \, d x= \int_{\partial \Omega} g (u-v) \, d \mathcal{H}^{N-1}, 
    \end{equation}
    for all $v\in W^{1,q}_{a}(\Omega)$.
\end{rem}

\begin{thm}\label{thm: existence} Under the assumptions of Proposition \ref{prop:behaviour},
    there exists at least one solution to problem \eqref{limiting}, in the sense of Definition \ref{defsol}.
\end{thm}
\begin{proof}
    Let $u_p$ be the solutions to problems \eqref{approximate} and $z\in L^{\infty}(\Omega)^N$ the vector field obtained from the previous Proposition. Then, we have that $\|z\|_{\infty}\le 1$. For any $\phi \in C_{0}^{\infty}(\Omega)$ by the weak formulation \eqref{weakform} we have that
\begin{equation*}
    \int_{\Omega} |\nabla u_p|^{p-2} \nabla u_p\cdot \nabla \phi  +a(x)|\nabla u_p|^{q-2} \nabla u_p  \cdot \nabla \phi \, d x =0\,.
\end{equation*}
Letting $p\to 1$ and using the previous Proposition we get that
\begin{equation}
\label{distributional}
    \int_{\Omega} z\cdot \nabla\phi +a(x)|\nabla u|^{q-2} \nabla u\cdot \nabla \phi \, d x =0 \,,\text{ for all }\phi\in C^\infty_0(\Omega)\,.
\end{equation}

\noindent Next, we prove that 
\begin{equation*}
z \cdot \nabla u=|\nabla u| \, \text{ a.e. in } \, \Omega.
\end{equation*}

\noindent Let $\phi \in C_{0}^{\infty}(\Omega)$ with $\phi\geq 0$ and take $u_p \phi$ as a test function  in \eqref{weakform}. Then
\begin{eqnarray*}
\int_{\Omega} \phi|\nabla u_p |^{p} d x+ \int_{\Omega} u_p |\nabla u_p |^{p-2} \nabla u_p \cdot \nabla \phi \, d x+\int_{\Omega} a(x)\phi|\nabla u_p |^q \, d x+\\
+\int_{\Omega} a(x) u_p|\nabla u_p |^{q-2} \nabla u_p \cdot \nabla \phi \, d x=0.
\end{eqnarray*}
By Young's inequality we have that
\[\int_\Omega\phi|\nabla u_p|\,d x\leq\frac{1}{p}\int_\Omega\phi |u_p|^pd x+\frac{p-1}{p}\int_\Omega \phi \,d x\]
and hence from the previous equation we get that
\begin{eqnarray}
\nonumber p\int_{\Omega} \phi|\nabla u_p | d x+ \int_{\Omega} u_p |\nabla u_p |^{p-2} \nabla u_p \cdot \nabla \phi \, d x
 +\int_{\Omega} a(x)\phi|\nabla u_p |^q \, d x+\\
 \label{final}+\int_{\Omega} a(x) u_p|\nabla u_p |^{q-2} \nabla u_p \cdot \nabla \phi \, d x \leq (p-1)\int_\Omega \phi \, d x\,.   
\end{eqnarray}
By Proposition \ref{prop:behaviour} we get that
\[u_p|\nabla u_p |^{p-2} \nabla u_p\stackrel{w}{\rightarrow} u z\text{ in }  L^{s}(\Omega)^N\,,\text{ for all }1\le s <\frac{N}{N-1}\]
\noindent and hence
\begin{equation*}
    \int_{\Omega} u_p \, (|\nabla u_p |^{p-2} \nabla u_p \cdot \nabla \phi )\, d x \to \int_{\Omega} u \, (z \cdot \nabla \phi )\, d x.
\end{equation*}
\noindent Next we have that
\begin{equation*}
 \int_{\Omega} u_p  a(x) |\nabla u_p |^{q-2} \nabla u_p \cdot \nabla \phi \, d x - \int_{\Omega} u a(x) |\nabla u |^{q-2} \nabla u  \cdot \nabla \phi \, d x =
 \end{equation*}
 \begin{eqnarray}
 \nonumber&=&\int_{\Omega} (u_p-u) a(x) |\nabla u_p |^{q-2} \nabla u_p  \cdot \nabla \phi \, d x+ \\
\label{second}&+&\int_{\Omega} u  a(x) \left(|\nabla u_p |^{q-2} \nabla u_p-|\nabla u |^{q-2} \nabla u \right)  \cdot \nabla \phi \, d x\,.
\end{eqnarray}
Again by Proposition \ref{prop:behaviour}
$$\nabla u_p\rightarrow\nabla u\text{ in }L^{q}_a(\Omega)^N\text{ and }|\nabla u_p |^{q-2} \nabla u_p\stackrel{w}{\rightarrow} |\nabla u |^{q-2} \nabla u\text{ in }L^{q'}_a(\Omega)^N $$
and thus the second summand in \eqref{second} converges to 0 as $p\rightarrow 1$.

For the first summand using H\"{o}lder's and Poincar\'{e}'s inequality \eqref{PoincareAq} we get that
\begin{eqnarray*}
    \int_{\Omega} (u_p-u)  a(x) |\nabla u_p |^{q-2} \nabla u_p \cdot \nabla \phi \, d x &\le& \|\nabla \phi\|_{\infty} \, \|u_p - u\|_{L^{q}_{a}} \, \|\nabla u_p\|_{L^{q}_{a}}^{\frac{q}{q'}}\\
    &\le& C\|\nabla u_p - \nabla u\|_{L^{q}_{a}} \, \|\nabla u_p\|_{L^{q}_{a}}^{\frac{q}{q'}} 
\end{eqnarray*}
and hence it also converges to 0.

\noindent Using the above and the weak lower semicontinuity  of the norm we pass to the limit in \eqref{final} and get that
\begin{eqnarray}
\nonumber\int_{\Omega} \phi|\nabla u |\, d x+ \int_{\Omega} uz \cdot \nabla \phi \, d x
 +\int_{\Omega} a(x)\phi|\nabla u |^q \, d x+\\
\label{final1}  +\int_{\Omega} a(x) u|\nabla u |^{q-2} \nabla u \cdot \nabla \phi \, d x \leq 0\,.   
\end{eqnarray}
If $v\in W^{1,q}_0(\Omega)$ then by density and using that $W^{1,q}_0(\Omega)\hookrightarrow W^{1,q}_{a,0}(\Omega)$, we get from \eqref{distributional} that  
\begin{equation}
\label{weakw1q}
    \int_{\Omega} z\cdot \nabla v +a(x)|\nabla u|^{q-2} \nabla u\cdot \nabla v \, d x =0\,.
\end{equation}
Since $u\phi\in W_{0}^{1,1}(\Omega)\cap W^{1,q}_{a, 0}(\Omega)$ we have by Proposition \ref{prop:MS} that there exists a sequence $(v_n)$ in $W^{1,q}_0(\Omega)$ such that 
\[v_n\to u\phi \text{ in } W^{1,1}(\Omega)\text{ and } \nabla v_n \to \nabla (u\phi)\text{ in } L^{q}_{a}(\Omega)^N\,.\]
Since each $v_n$ satisfies \eqref{weakw1q}, passing to the limit we get that
\[\int_{\Omega} z \cdot \nabla (u\phi) \, d x+ \int_{\Omega} a(x) |\nabla u|^{q-2} \nabla u \cdot \nabla (u\phi) \, d x = 0\,.\]

\noindent Combining this with \eqref{final1} we get that
\begin{equation*}
\int_{\Omega} \phi  \, |\nabla u | \, d x \le  \int_{\Omega} \phi  \, z \cdot \nabla u  \, d x, \, \text{ for all }\, \phi \in C^{\infty}_{0}(\Omega)
\end{equation*}
and hence
$$|\nabla u|\le z\cdot \nabla u\text{ a.e. in } \Omega\,.$$ 
Since on the other hand $\|z\|_{\infty}\le 1$ implies that $z \cdot \nabla u \le |\nabla u|$ we infer that
\begin{equation*}
    z \cdot \nabla u = |\nabla u| \, \text{ a.e. in }\, \Omega.
\end{equation*}
To conclude the proof, we need to show that 
\begin{equation*}
    \gamma(\zeta) = g \,\,\,\, \mathcal{H}^{N-1} \, \text{- a.e. on } \, \partial \Omega. 
\end{equation*}
We will use the space of test functions $D(\overline{\Omega})=\{ v|_{\Omega} : v \in D(\mathbb{R}^N)\}$. Let $v\in D(\overline{\Omega}) \hookrightarrow W^{1,\theta_p}(\Omega)$ as a test function in \eqref{weakform}
\begin{equation}
\nonumber
\int_{\Omega} \left(|\nabla u_p|^{p-2} \nabla u_p +a(x)|\nabla u_p|^{q-2}\nabla u_p \right)  \nabla v \, dx = \int_{\partial \Omega} g v \, d\mathcal{H}^{N-1}.
\end{equation}
Passing to the limit as $p\to 1$, by Proposition \ref{prop:behaviour} we get
\begin{equation}
\nonumber
\int_{\Omega} \left(z +a(x)|\nabla u|^{q-2}\nabla u \right)  \nabla v \, dx = \int_{\partial \Omega} g v \, d\mathcal{H}^{N-1},
\end{equation}
and applying Gauss-Green formula \eqref{GaussGreen} we get
\begin{equation}
\nonumber
\int_{\Omega } v \,\mathrm{div}\left(z +a(x)|\nabla u|^{q-2}\nabla u \right) \,  dx + \langle  \gamma (\zeta),  v \rangle_{W^{-\frac{1}{q},q} \times W^{1-\frac{1}{q},q}} = \int_{\partial \Omega} g v \, d\mathcal{H}^{N-1},
\end{equation}
where $\zeta =z+a(x)|\nabla u|^{q-2} \nabla u $. Using the fact that $\mathrm{div}\left(z +a(x)|\nabla u|^{q-2}\nabla u \right)=0$ we get
\begin{equation}
\nonumber
\langle  \gamma (\zeta),  v \rangle_{W^{-\frac{1}{q},q} \times W^{1-\frac{1}{q},q}} = \int_{\partial \Omega} g v \, d\mathcal{H}^{N-1} \,\, \text{ for all } \, v\in D(\overline{\Omega})
\end{equation}
or  
\begin{equation}
\nonumber
\langle  \gamma (\zeta)-g,  v \rangle_{W^{-\frac{1}{q},q} \times W^{1-\frac{1}{q},q}} = 0 \,\, \text{ for all } \, v\in D(\overline{\Omega}),
\end{equation}
since $L^{\infty}(\partial \Omega) \hookrightarrow    W^{-\frac{1}{q},q} (\partial \Omega)$. By density of $D(\overline{\Omega})$ in the space $W^{1,q}(\Omega)$, we get that
\begin{equation*}
    \gamma(\zeta)=g \,\,\,\, \mathcal{H}^{N-1} \, \text{- a.e. on } \, \partial \Omega. 
\end{equation*}
\end{proof}

\begin{rem}\label{solutions}
    Using the weak formulation \eqref{weaklimiting}, it is easy to see that if $u_1,u_2$ are two weak solutions of problem \eqref{limiting} then they differ up to a constant. To this end, let
    $u_1,u_2\in W^{1,q}_{a}(\Omega)$ be two solutions and $z_1,z_2 \in L^{\infty}(\Omega)^N$ the corresponding vector fields. We can use each one of $u_1$ and $u_2$ as a test function in the weak formulation \eqref{weaklimiting} for other one and add the two equations to get
    \begin{eqnarray*}
\int_\Omega |\nabla u_1| \, d x&+&\int_\Omega |\nabla u_2| \, d x -\int_\Omega z_1 \cdot \nabla u_2 \, d x -\int_\Omega z_2 \cdot \nabla u_1 \, d x +\\
&+ &\int_\Omega a(x) \left(|\nabla u_1|^{q-2} \nabla u_1-|\nabla u_2|^{q-2} \nabla u_2 \right) \nabla (u_1-u_2)\, d x =0.\
\end{eqnarray*}
    Taking into account that 
    \begin{equation*}
        \int_{\Omega} z_2 \cdot \nabla u_1 \, d x \le \int_{\Omega} |\nabla u_1| \, d x 
    \end{equation*}
    and 
     \begin{equation*}
        \int_{\Omega} z_1 \cdot \nabla u_2 \, d x \le \int_{\Omega} |\nabla u_2| \, d x 
    \end{equation*}
we get
\begin{equation*}
\int_\Omega a(x) \left(|\nabla u_1|^{q-2} \nabla u_1-|\nabla u_2|^{q-2} \nabla u_2 \right) \cdot \nabla (u_1-u_2)\, d x \le 0\,.
\end{equation*}
and since the integrand is non-negative this implies that
\begin{equation*}
\int_\Omega a(x) \left(|\nabla u_1|^{q-2} \nabla u_1-|\nabla u_2|^{q-2} \nabla u_2 \right) \cdot \nabla (u_1-u_2)\, d x = 0\,.
\end{equation*}
Hence we conclude that
\[\nabla u_1=\nabla u_2\,,\text{ a.e. in }\Omega,\]
which implies that $u_1=u_2+c$ a.e in $\Omega$, for some $c\in \mathbb{R}$. Note that as the solution in Theorem \ref{thm: existence} is obtained as the limit $u_p \to u$ in $L^1(\Omega)$ and $\int_{\Omega}u_p \, dx=0$ we have that $\int_{\Omega}u \, dx=0$. Hence the solution is unique on the space $V(\Omega)$.
\end{rem}

The last result of this paper is to obtain a variational characterization of the solutions to the problem \eqref{limiting} through the minimization
 problem for the functional

\begin{eqnarray}\label{min}
\nonumber
   \mathcal{I}(u)= \int_{\Omega} |\nabla u| \, dx + \frac{1}{q} \int_{\Omega} a(x) |\nabla u|^q \, dx - \int_{\partial \Omega} g u \, d \mathcal{H}^{N-1},
\end{eqnarray}
on the set
\begin{equation*}
    V(\Omega)=\bigg\{ u \in  W^{1,q}_{a}(\Omega): \int_{\Omega} u \, dx=0 \bigg\}.
\end{equation*}

\begin{prop}
    The function $u\in V(\Omega)$ is the unique weak solution to problem \eqref{limiting} if and only if is the unique minimizer of $\mathcal{I}$.
\end{prop}

\begin{proof}
    By Theorem \ref{thm: existence} and Remark \ref{solutions} there exists a solution $u$, which is unique on the set $V(\Omega)$. Moreover, the functional $\mathcal{I}$ is strictly convex on the set $V(\Omega)$ hence, it is enough to show that $u$ is a minimizer of $\mathcal{I}$. To this end, let $z\in L^{\infty}(\Omega)^N$ be the vector field given by Theorem \ref{thm: existence} and $v\in V(\Omega)$. By the weak formulation \eqref{weakformlim}, we have that
    \begin{equation*}
         \int_{\Omega} |\nabla u| \, d x-\int_{\Omega} z \cdot \nabla v \, d x + \int_{\Omega} a(x) |\nabla u|^{q-2} \nabla u \cdot \nabla (u-v) \, d x= \int_{\partial \Omega} g (u-v) \, d \mathcal{H}^{N-1}.
    \end{equation*}
    Using Young's inequality and the fact that $\|z\|_{\infty}\le 1$ we get
   \begin{eqnarray*}
\int_{\Omega} |\nabla u| \, d x &+&\int_{\Omega} a(x) |\nabla u|^q\, d x - \int_{\partial \Omega} g u \, d \mathcal{H}^{N-1}\\&= &\int_{\Omega} z \cdot \nabla v \, d x+ \int_{\Omega} a(x) |\nabla u|^{q-2} \nabla u \cdot \nabla v \, d x- \int_{\partial \Omega} g v \, d \mathcal{H}^{N-1}\\
&\leq & \int_{\Omega} |\nabla v| \, d x+\frac{1}{q'} \int_{\Omega} a(x) |\nabla u |^{q} d x + \frac{1}{q} \int_{\Omega} a(x) |\nabla v |^{q} d x\\
&-& \int_{\partial \Omega} g v \, d \mathcal{H}^{N-1}\,,
\end{eqnarray*}
  which leads to 
  \begin{equation*}
      \mathcal{I}(u)\le \mathcal{I}(v),
  \end{equation*}
  for all $v\in V(\Omega)$ and the proof is completed.
\end{proof}

\subsection*{Acknowledgment}
We would like to thank Prof. N. S. Papageorgiou for the valuable and fruitful discussions. We would also like to thank the anonymous referee for his/hers comments and suggestions that helped improve both the content and the presentation of this paper.

\end{document}